\documentclass[12pt,leqno]{amsart}
\usepackage{amscd,times,fullpage}
\usepackage{ucs}
\usepackage[utf8x]{inputenc}
\usepackage{amsmath}
\usepackage{extarrows}
\usepackage{amsthm,amsmath}
\usepackage{amssymb}
\usepackage{graphicx}
\usepackage{amsthm}
\usepackage{enumerate}

\usepackage[all]{xy}
\usepackage{tikz}
\usetikzlibrary{positioning} 
\usetikzlibrary{arrows,decorations.pathmorphing,shapes} 
\usepackage{tikz-cd}
\usetikzlibrary{cd}
\usetikzlibrary{matrix}

\newcommand{\lra}{\longrightarrow}             
\newcommand{\R}{\mathbb{R}} 
\newcommand{\B}{\text{B}}                   
\newcommand{\Z}{\mathbb{Z}}
\newcommand{\q}{\mathbb{Q}}
\newcommand{\CP}{\mathbb{C}\mathrm{P}}
\newcommand{\N}{\mathbb{N}}

\newcommand{\C}{\mathbb{C}}

\newcommand{\di}{{\operatorname{d}}}
\newcommand{\betti}{{\operatorname{b}}}

\newcommand{\Id}{\operatorname{Id}}

\newcommand{\coker}{{\operatorname{coker}}}
\newcommand{\F}{\mathcal{F}}

\newcommand{\Pg}{\mathcal{P}}
\newcommand{\Kg}{\mathcal{K}}
\newcommand{\Sg}{\mathcal{S}}
\newcommand{\D}{\mathcal{D}}
\newcommand{\TM}{\text{T} M}

\newtheorem{theor}{Theorem}
\newtheorem{prop}[theor]{Proposition}

\newtheorem{lem}[theor]{Lemma}
\newtheorem{cor}[theor]{Corollary}

\newtheorem{ex}[theor]{Example}
\newtheorem{remark}[theor]{Remark}

\begin{document}

\title{Sasaki versus K\"ahler groups}

\author{D.~Kotschick}
\address{Mathematisches Institut, Ludwig-Maximilians-Universit\"at
M\"unchen, Theresienstr.~39, 80333 M\"unchen, Germany}
\email{dieter@math.lmu.de}

\author{G.~Placini}
\address{Dipartimento di Matematica e Informatica “Ulisse Dini”, Universit\`{a} di Firenze, viale Morgagni 67/A, 50134 Firenze, Italy}
\email{giovanni.placini@unifi.it, giovanni.placini@unica.it}

\date{\today ; {\copyright \ D.~Kotschick and G.~Placini 2022}}




\begin{abstract}
We study fundamental groups of compact Sasaki manifolds and show that compared 
to K\"ahler groups, they exhibit rather different behaviour. This class of groups is
not closed under taking direct products, and there is often an upper bound
on the dimension of a Sasaki manifold realising a given group. 
The richest class of Sasaki groups arises in dimension $5$.
\end{abstract}
 
\maketitle

\section{Introduction}

Sasaki geometry is often considered as the odd-dimensional analogue of K\"ahler geometry. In this paper 
we examine this analogy with regard to the fundamental group. It is well known now that K\"ahler groups,
 that is the fundamental groups of compact K\"ahler manifolds, form a very special subclass of of the 
 class of all finitely presentable groups; cf.~\cite{ABCKT}. 
 In their foundational monograph~\cite{Boyer} on Sasaki geometry, Boyer and Galicki suggested that 
 one study Sasaki groups, that is, fundamental groups of compact Sasaki manifolds, systematically. 
 They noted that groups with virtually odd first Betti number cannot be Sasaki, and then suggested that perhaps 
 the analogy with K\"ahler groups might stop there, cf.~\cite[p.~236]{Boyer}. 
 Nevertheless, in the papers written on Sasaki groups since then, the point of view taken was usually 
 that of pushing the analogy further, by proving that certain constraints known for K\"ahler groups also apply 
 to Sasaki groups; see Chen~\cite{Chen}, Biswas et.~al.~\cite{BFMT,BM} and Kasuya~\cite{Kas2,Kas1}.
 
 In this paper we will also prove some results which extend the analogy between
 the K\"ahler and Sasaki cases. However, our main point is that Sasaki groups actually behave 
 rather differently from K\"ahler groups. The reason that many results about K\"ahler 
 groups extend to Sasaki groups is that within certain restricted classes of groups all
 Sasaki groups are in fact K\"ahler, and, therefore, within those classes of groups one does 
 not encounter the important differences between the K\"ahler and Sasaki situations.

 \subsection{K\"ahler and projective groups}
 
 Let us denote by $\Kg_n$ the class of fundamental groups of closed K\"ahler manifolds of complex 
 dimension $n$, and by $\Pg_n\subset\Kg_n$ the subclass of groups realised by smooth complex projective 
algebraic varieties. By taking products with $\C P^1$ one sees that there are 
natural inclusions $\Kg_n\subset\Kg_{n+1}$ and $\Pg_n\subset\Pg_{n+1}$. The latter inclusion is in 
fact an equality for $n\geq 2$ because of the Lefschetz hyperplane theorem. We therefore denote 
$\Pg_{n\geq 2}$ by $\Pg$. It is an open question whether the inclusion $\Kg_n\subset\Kg_{n+1}$
is ever strict for $n>1$.

For $n=1$ we have $\Kg_1=\Pg_1$, since every compact Riemann surface is an algebraic curve. These
fundamental groups are just the orientable surface groups.
A classical result of Kodaira says that $\Kg_2=\Pg$, and a recent result of Claudon, H\"oring and Lin~\cite{CHL}
gives $\Kg_3=\Pg$. Therefore, $\Kg_2=\Kg_3$, a statement for which no direct proof is known, other 
than arguing that both sides equal $\Pg$.

We can summarise all this information in the diagram of inclusions depicted in Figure~\ref{Kaehler}.

\begin{figure} \begin{tikzcd} [column sep={3em,between origins},
  row sep={3em,between origins},]
\Kg_1\arrow[d, phantom , sloped , "="] & \subsetneq & \Kg_2\arrow[d, phantom , sloped , "="] & = & \Kg_3\arrow[d, phantom , sloped , "="] & \subset & \Kg_4 & \subset & \Kg_5 & \subset & \ldots \\
\Pg_1 & \subsetneq & \Pg & = & \Pg & = & \Pg\arrow[u, phantom, sloped , "\subset"] & = & \Pg\arrow[u, phantom, sloped , "\subset"]   & = & \ldots
\end{tikzcd}
\caption{K\"ahler and projective groups}\label{Kaehler}
\end{figure}

The class of K\"ahler groups is the union of all the $\Kg_n$. There is a strict horizontal inclusion
$\Kg_n\subsetneq\Kg_{n+1}$ for some $n\geq 3$ if and only if one\footnote{and, therefore, infinitely many} 
of the vertical inclusions is strict,
meaning that there would be more K\"ahler groups than projective ones. If this were the case, then the 
smallest dimension in which a K\"ahler group could arise as the fundamental group of a closed 
K\"ahler manifold would be an interesting invariant of (non-projective) K\"ahler groups. Claudon, H\"oring and 
Lin\footnote{and perhaps others as well} have 
conjectured that all K\"ahler groups are projective, cf.~\cite[Conjecture~1.1]{CHL}, which would mean
that there are no vertical strict inclusions in this diagram.

Finally note that the classes of K\"ahler or projective groups are both closed under taking subgroups 
of finite index, and that they are also closed under direct products.

 \subsection{Sasaki and projective groups}
 
 Let us denote by $\Sg_n$ the class of fundamental groups of closed Sasaki manifolds of 
 dimension $n$. Taking coverings shows that this is closed under passing to subgroups of 
 finite index. However, there is no obvious way to realise direct products, since the Cartesian
 product of two Sasaki manifolds is of even dimension, and therefore certainly not Sasaki.
 This problem not only obstructs the realisation of direct products of groups, it also 
 means that there is no straightforward way to get an inclusion of $\Sg_n$ into $\Sg_m$ for 
 a larger $m$, be it $n+2$ or some larger value.
 
 It is known from work of Kasuya~\cite{Kas1} that these issues cannot be resolved when $n=3$.
 The Sasaki $3$-manifolds with infinite fundamental groups are essentially the circle bundles 
 with non-zero Euler classes over surfaces of positive genus, cf.~Geiges~\cite{Geiges}.
 They are not $1$-formal, and, in fact, have non-trivial Massey triple products defined
 on $H^1$. Now Kasuya~\cite[Theorem~1.1]{Kas1} proved that Sasaki manifolds of dimension $\geq 5$
 are $1$-formal, equivalently the Malcev algebra of a Sasaki group is quadratically presented.
 Therefore no infinite group in $\Sg_3$ will occur in $\Sg_n$ for any $n\geq 5$.
The same conclusion holds for any direct product of groups that has an infinite factor from $\Sg_3$.

Kasuya's result suggests that groups in $\Sg_3$ are completely unlike high-dimensional Sasaki groups,
and one might hope that once one discards the $3$-dimensional case, things will be more or 
less uniform, as expected by analogy with the K\"ahler case. Of course, Kasuya's result 
about $1$-formality further strengthens the analogy between K\"ahler groups and high-dimensional 
Sasaki groups. 

Our first result is that in all odd dimensions there are Sasaki groups which do not occur in any 
larger dimension, so that the $3$-dimensional case is actually not special in this regard.
\begin{theor}\label{t1}
Every $\Sg_{n}$ contains groups which are not contained in any $\Sg_{m}$ with $m>n$.
\end{theor}
This means that the problem caused by the absence of products in the Sasaki category
cannot be resolved, and there is no inclusion between the $\Sg_n$ going up in dimension. However,
we can adapt the Lefschetz hyperplane theorem to obtain an inclusion going down in dimension:
\begin{theor}\label{t2}
There is an inclusion $\Sg_{n}\subset\Sg_{n-2}$ for all $n\geq 7$. 
\end{theor}
By the previous theorem, all these inclusions are strict. The combination of the two theorems shows 
that $\Sg_5$ is the union of all $\Sg_{n\geq 5}$, and that it is strictly larger than any union taken by leaving 
out $\Sg_5$. It is somewhat counterintuitive that the largest and most interesting class of fundamental
groups occurs in a rather small dimension, in fact the smallest dimension beyond the exceptional $n=3$.

\begin{figure} \begin{tikzcd} [column sep={3em,between origins},
  row sep={3em,between origins},]
\Sg_3 &  & \Sg_5 & \supsetneq & \Sg_7 & \supsetneq & \Sg_9 & \supsetneq & \Sg_{11} & \supsetneq & \ldots \\
 & & \Pg\arrow[u, phantom, sloped , "\subsetneq"] & = & \Pg\arrow[u, phantom, sloped , "\subsetneq"] & = & \Pg\arrow[u, phantom, sloped , "\subsetneq"] & = & \Pg\arrow[u, phantom, sloped , "\subsetneq"]   & = & \ldots
\end{tikzcd}
\caption{Sasaki and projective groups}\label{Sasaki}
\end{figure}
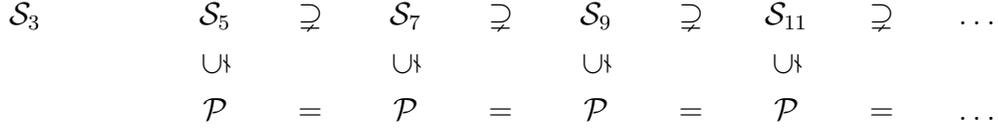

We also show that projective groups are Sasaki in all dimensions $>3$.
\begin{theor}\label{t3}
The class $\Sg_{n}$ of Sasaki groups of dimension $n$ contains all projective 
groups for every $n\geq 5$. 
\end{theor}
Without control on the dimension, this was proved previously
by Chen~\cite[Proposition~1.2.]{Chen}. His argument does not yield the most interesting case of $\Sg_5$.
Note that by Theorem~\ref{t1} all the inclusions $\Pg\subset\Sg_{n}$ in Theorem~\ref{t3} are strict.
For Sasaki groups the analog of the diagram in Figure~\ref{Kaehler} is the diagram in Figure~\ref{Sasaki}, 
summarizing our discussion so far, and leaving $\Sg_3$ out in left field in splendid isolation.
In light of this diagram it seems interesting to investigate the intersection of all the $\Sg_n$ for $n\geq 5$.
This intersection contains all the projective groups, and one might be tempted to conjecture that it equals $\Pg$. 
However, Theorem~\ref{t3} is true not only for the fundamental groups of smooth complex projective algebraic 
varieties, but also for the potentially larger class of orbifold fundamental groups of cyclic polarised projective 
orbifolds in the sense of Ross and Thomas~\cite[Definition~2.7.]{RT}, see Remark~\ref{orbifoldremark} below. 
Therefore, at least these groups are contained in all $\Sg_n$ for $n\geq 5$.

\subsection{Extending restrictions on K\"ahler groups to Sasaki groups}

By Theorem~\ref{t3} all projective groups are Sasaki, and, therefore, if all K\"ahler groups are projective, then
all K\"ahler groups are Sasaki. However, there are many Sasaki groups which are not K\"ahler, for 
example all the ones constructed in the proof of Theorem~\ref{t1} have this property, 
see Section~\ref{s:SasakiNotKaehler} . Nevertheless, within certain restricted classes of groups 
it turns out that all Sasaki groups are in fact projective. In this direction we will prove the 
following:
\begin{theor}\label{t4}
In the following classes of groups all Sasaki groups in $\Sg_5$ are in fact projective:
\begin{enumerate}
\item \label{nocentre} torsion-free groups with trivial centre,
\item \label{hyperbolic} torsion-free hyperbolic groups,
\item \label{Schreier}  torsion-free Schreier groups,
\item \label{rankone}  fundamental groups of non-positively curved closed manifolds of rank one,
\item \label{SSNT}  fundamental groups of compact locally symmetric spaces of non-compact type,
\item \label{3D}  fundamental groups of three-manifolds.
\end{enumerate}
\end{theor}
This applies in particular to fundamental groups of manifolds of constant negative curvature,
see Corollary~\ref{c:hyperbolic}.
The philosophy behind Theorem~\ref{t4} also explains other results that have been proved about 
Sasaki groups. The reason certain restrictions on K\"ahler groups hold for Sasaki groups is that
Sasaki groups satisfying those restrictions are in fact projective. For example, {\it a posteriori}, the recent 
result of Biswas and Mj~\cite{BM} can be paraphrased as saying that Sasaki groups of deficiency
at least $2$ must be projective, compare~\cite{Kot2}. For torsion-free groups this is an immediate consequence 
of Statement~\ref{Schreier} in Theorem~\ref{t4}.

There are other restrictions on K\"ahler groups that can be extended to the Sasaki case, but 
that do not fit neatly into this philosophy. A case in point is the 
following Sasaki analog of the well known theorem of Johnson and Rees~\cite{johnsonrees87} about K\"ahler groups. 
\begin{theor}\label{t5}
Let $\Gamma_1$ and $\Gamma_2$ be two groups. Assume $f_i\colon \Gamma_i\lra Q_i$ are non-trivial quotients with $\vert Q_i\vert=m_i<\infty$ for both $i=1,2$. Then the following two statements hold:
\begin{enumerate}[(a)]
\item $\Gamma_1*\Gamma_2$ is not Sasaki, \label{johnsonparta}
\item moreover, for any group $H$ the product $(\Gamma_1*\Gamma_2)\times H$ is not Sasaki. \label{johnsonpartb}
\end{enumerate} 
\end{theor}

\subsection{Structure of the paper}
In Section~\ref{sectiontop} we introduce notation and recall some facts on the topology of Sasaki manifolds 
that are crucial for our arguments.
Section~\ref{s:Sasakiprojective} is dedicated to the relationship between Sasaki and projective groups, in particular
it contains the proofs of Theorems~\ref{t2} and \ref{t3}. We also explain how the proof of Theorem~\ref{t3} adapts
to $K$-contact manifolds.
In Section~\ref{s:SasakiNotKaehler} we discuss Sasaki groups which are not projective, and not even K\"ahler,
and which can only be realised as Sasaki groups on manifolds with explicit dimension bounds. This in particular
proves Theorem~\ref{t1}.
Finally Section~\ref{s:further} contains the proofs of Theorems~\ref{t4} and \ref{t5}, and some corollaries and 
applications.

\subsection*{Acknowledgments} 
The first named author would like to thank B.~Claudon and P.~Py for useful discussions and communications.

\section{Sasaki manifolds and the associated group extensions}\label{sectiontop}

We begin with some definitions and known results; 
for a more comprehensive treatment we refer to the monograph by Boyer and Galicki~\cite{Boyer}. 
Unless otherwise stated, all manifolds are assumed to be smooth, closed and oriented.

A \textit{K-contact structure} $(M,\eta,\phi)$ on a manifold $M$ consists of a contact form $\eta$ and an endomorphism 
$\phi$ of the tangent bundle 
$\TM$ satisfying the following properties:
\begin{enumerate}
\item[$\bullet$] $\phi^2=-\Id+R\otimes\eta$ where $R$ is the Reeb vector field of $\eta$,
\item[$\bullet$] $\phi_{\vert\D}$ is an almost complex structure compatible with the symplectic form $\di\eta$ on $\D=\ker\eta$,
\item[$\bullet$] the Reeb vector field $R$ is Killing with respect to the metric $g(\cdot,\cdot)=\di\eta(\phi\cdot,\cdot)+\eta(\cdot)\eta(\cdot)$.
\end{enumerate}   
Given such a structure one can consider the almost complex structure $I$ on the Riemannian cone $\big( M\times(0,\infty),t^2g+\di t^2\big)$ given by
\begin{enumerate}
\item[$\bullet$] $I=\phi$ on $\D=\ker\eta$,
\item[$\bullet$] $I(R)=t\partial_t$.
\end{enumerate} 
A \textit{Sasaki structure} is a K-contact structure $(M,\eta,\phi)$ such that the associated almost complex structure $I$ on the Riemannian cone is integrable.

An important example of a Sasaki (resp.~K-contact) structure is obtained by the \textit{Boothby--Wang fibration} 
$M$ over a K\"ahler (resp.~almost-K\"ahler) manifold $(X,\omega)$ with $\omega$ representing an integral class~\cite{BW}, 
that is, the principal $S^1$-bundle $\pi\colon M\longrightarrow X$ with Euler class $[\omega]$ and connection 
$1$-form $\eta$ such that $\pi^*(\omega)=\di\eta$.

A contact form is called \textit{regular} (respectively \textit{quasi-regular}, \textit{irregular}) 
if its Reeb foliation is such. 
Rukimbira~\cite{Ruk} proved that any irregular Sasaki structure can be deformed to a quasi-regular one. 
Moreover, the geometry of quasi-regular structures is described by the following 
(cf.~\cite[Theorem 7.5.1]{Boyer}):
\begin{theor}[Structure Theorem]\label{structure}
Let $(M,\eta,\phi)$ be a quasi-regular Sasaki structure and let $\vert X\vert$ be the space of leaves of the Reeb folation. 
Then $\vert X\vert$ carries the structure of a projective orbifold $X$ with an integral K\"ahler class 
$[\omega]\in H^2_{orb}(X;\Z)$, and $\pi\colon M\longrightarrow X$ is the principal $S^1$-orbibundle with 
connection $1$-form $\eta$ such that $\pi^*\omega=\di\eta$. 
Moreover, if $\eta$ is regular then $\pi\colon M\longrightarrow X$ is a principal 
$S^1$-bundle over a smooth projective manifold.
\end{theor}

The orbifold homotopy, homology and cohomology groups are defined using Haefliger's classifying 
space~\cite{Haefliger}, cf.~\cite[Section 4.3]{Boyer}.

Since the leaf holonomy groups of the Reeb foliation are always cyclic, the orbifold $X$ is a 
cyclic orbifold. In particular, it is normal and $\q$-factorial. Moreover, because $M$ is a 
smooth manifold, the cohomology class of $\omega$ being integral means that $X$ is a polarised orbifold 
 in the sense of Ross and Thomas~\cite[Definition~2.7]{RT}.
Throughout this paper, the only orbifolds we consider are these cyclic polarised orbifolds.
Starting from one of these orbifolds, one has the following converse to the structure theorem
(cf.~\cite[Theorem 7.5.2]{Boyer}), arising from the orbifold version of the Boothby--Wang 
construction.
\begin{theor}\label{conversestructure}
Let $X$ be a cyclic projective orbifold equipped with a polarisation defined by an integral K\"ahler class 
$[\omega]\in H^2_{orb}(X;\Z)$. Then the total space of the principal $S^1$-orbibundle 
$\pi\colon M\longrightarrow X$ with Euler class $[\omega]$
is a manifold that can be equipped with a quasi-regular Sasaki structure such that $\pi^*\omega=\di\eta$,
where $\eta$ is the contact form.
\end{theor}
In this particular case the total space of the orbibundle is actually smooth, because we started with a polarised cyclic 
orbifold in the sense of Ross and 
Thomas~\cite[Definition~2.7]{RT}, which means in particular that all local uniformizing groups inject into $S^1$,
so that smoothness of the total space follows from~\cite[Lemma~4.2.8]{Boyer}.

The Structure Theorem~\ref{structure} together with Rukimbira's result~\cite{Ruk} implies that the fundamental group of a Sasaki manifold $M$ 
always fits into the long exact sequence of homotopy groups associated to the principal $S^1$-orbibundle $\pi\colon M\rightarrow X$, namely:
$$
\ldots\longrightarrow\pi_2^{orb}(X)\stackrel{\partial}{\longrightarrow}\pi_1(S^1)\longrightarrow\pi_1(M)\longrightarrow\pi_1^{orb}(X)\longrightarrow 0 \ .
$$
Thus $\pi_1(M)$ fits into a short exact sequence of the following type:
\begin{equation}\label{ses}
0\longrightarrow C=\coker \partial\longrightarrow\pi_1(M)\longrightarrow\pi_1^{orb}(X)\longrightarrow 0 \ .
\end{equation} 
The group $C$ is a -- possibly trivial -- cyclic subgroup of the centre of $\pi_1(M)$  determined by the Euler class $e=[\omega]\in H^2_{orb}(X;\Z)$ since the map $\partial$ 
factors as 
\begin{equation}\label{EqMapFact}
\begin{tikzcd}
\pi_2^{orb}(X) \arrow[r, "\partial"] \arrow[dr, "\psi"'] & \Z\cong\pi_1(S^1) \\
& H_2^{orb}(X;\Z ) \arrow[u, "\langle - ;e\rangle"'] 
\end{tikzcd}
\end{equation}
where $\langle - ;e\rangle$ is the evaluation of the Euler class and $\psi$ is the Hurewicz homomorphism; see for instance~\cite{Ham}.

\begin{remark}\label{euler}
	Let $\B\pi_1^{orb}(X)$ be the classifying space constructed from the orbifold classifying space  $\B X$ by attaching cells of 
	dimension larger than $2$. The inclusion map $\iota\colon\B X\lra \B\pi_1^{orb}(X)$ induces an isomorphism 
	$\iota^*\colon H^1(\B \pi_1^{orb}(X);\Z)\lra H^1_{orb}(X;\Z)$  and an injective homomorphism
	$\iota^*\colon H^2(\B \pi_1^{orb}(X);\Z)\hookrightarrow H^2_{orb}(X;\Z)$.

	In the case where $C\cong \Z$ in~\eqref{ses} one can identify the Euler class $e\in H^2_{orb}(X;\Z)$ of the principal $S^1$-orbibundle 
	$p\colon M\longrightarrow X$ with the characteristic class $c\in H^2(\pi_1^{orb}(X);\Z)$ of the central extension via the pullback of the inclusion $\iota\colon\B X\longrightarrow \B\pi_1^{orb}(X)$. 
\end{remark}

Let us now give some more details about the extension~\eqref{ses}.
In particular, we want to relate the orbifold fundamental group $\pi_1^{orb}(X)$ to a genuine projective group. 
Note that the map $p\colon BX\lra X$ from the orbifold classifying space to the underlying topological space 
induces a surjective map $p_*$ at the level of fundamental groups. 
Moreover, the kernel of $p_*$ is normally generated by loops around the irreducible divisors $D$ contained in the singular set of $X$.
These loops represent torsion elements of order $m$, the ramification index of $D$.
Therefore, the kernel $K$ of the map $p_*\colon\pi_1^{orb}(X)\lra \pi_1(X)$ is generated by (possibly infinitely many) torsion elements.
Now the cohomology of $K$ with coefficients in $\R$ is trivial in all positive degrees because $\R$ is $m$-divisible for all $m$.
Thus the Lyndon--Hochschild--Serre spectral sequence of $K\subset \pi_1^{orb}(X)$ yields an isomorphism
$H^*(\pi_1^{orb}(X);\R)\cong H^*(\pi_1(X);\R)$.
Moreover, $X$ admits a resolution of singularities which does not change the fundamental group by a result of Koll\'ar, see \cite[Theorem~7.5.2]{kollar93}.
Thus the real cohomology ring of $\pi_1^{orb}(X)$ is that of a projective group.
Notice that whenever $C\neq\Z$ we have an isomorphism $H^*(\pi_1^{orb}(X);\R)\cong H^*(\pi_1(M);\R)$. In this instance $\pi_1(M)$ itself has the real cohomology ring of the projective group $\pi_1(X)$.
We summarize this discussion in a lemma for future reference.
\begin{lem}\label{LemCohomSasakiGroup}
	For any quasi-regular structure $\pi\colon M\lra X$ on a Sasaki manifold $M$ one has the diagram
	\begin{equation}\label{EqDiagramSasakiGroups}
		\begin{tikzcd}
			&             &                           & K \arrow[d, hook]&\\
			0 \arrow[r] & C \arrow[r] & \Gamma \arrow[r, "\pi_*"] & \pi_1^{orb}(X) \arrow[r] \arrow[d, two heads, "p_*"] & 0  \\
			&             &                           & \pi_1(X) &    
		\end{tikzcd} 
	\end{equation}
	where $\Gamma=\pi_1(M)$ and $C$ is cyclic.
	Moreover, $\pi_1(X)$ is a projective group and the kernel $K$ of $p_*$ is generated by torsion elements so that $H^*(\pi_1^{orb}(X);\R)\cong H^*(\pi_1(X);\R)$.
	If in addition $C\neq\Z$, then $H^*(\pi_1(X);\R)\cong H^*(\Gamma;\R)$.
\end{lem}

We will also need the following.
\begin{lem}\label{LemNonDegPairingOrbi}
	Let $\pi\colon M\lra X $ be the principal orbibundle associated to a quasi-regular Sasaki structure.
	Then there is a non-degenerate skew-symmetric bilinear pairing
	\begin{align}\label{EqNonDegPairingOrbi}
		H^1(\pi_1^{orb}(X );\R)\times H^1(\pi_1^{orb}(X );\R)\lra\R
	\end{align}
	which factors through the cup product
	\begin{align*}
		H^1(\pi_1^{orb}(X );\R)\times H^1(\pi_1^{orb}(X );\R)\xlongrightarrow{\cup} H^2(\pi_1^{orb}(X );\R).
	\end{align*}
\end{lem}
\begin{proof}
	In the quasi-regular case the basic cohomology ring $H_B^*(\F;\R)$ of the Reeb fibration, see~\cite[Section~7.2]{Boyer}, 
	coincides with the orbifold cohomology ring with real coefficients $H_{orb}^*(X ;\R)$.
	Moreover, the Hard Lefschetz Theorem holds for the basic cohomology of a Sasaki manifold, 
	see~\cite{elkacimialaoui90} and~\cite[Theorem~7.2.9]{Boyer}.
	The claim then follows by composing the non-degenerate bilinear map coming from the transverse Hard Lefschetz Theorem 
	with the map $\iota^*\colon H^*(B\pi_1^{orb}(X))\lra H^*_{orb}(X)$ given in Remark~\ref{euler}.
\end{proof}
\begin{remark}
This lemma is crucial for the proof of Theorem~\ref{t5}.
	In~\cite{cappellettimontanodenicolayudin15} a different version of the Hard Lefschetz Theorem is proved for Sasaki 
	manifolds. However, the non-degenerate bilinear pairing constructed in~\cite{cappellettimontanodenicolayudin15}  
	does not factor through the cup product in group cohomology, and is therefore not suitable for our purposes.
	
\end{remark}

\section{Sasaki groups from projective groups}\label{s:Sasakiprojective}

In this section we  prove Theorems~\ref{t2} and~\ref{t3}, and give some variations on the latter.

\subsection{The Lefschetz hyperplane theorem for Sasaki manifolds}

In algebraic geometry, the Lefschetz hyperplane theorem is a statement relating the homotopy
groups of a complex projective variety to those of a generic hyperplane section. In the Sasaki
context we consider quasi-regular structures defining orbifold bundles over projective orbifolds,
and then take hyperplane sections of the base orbifold and restrict the orbifold bundle to such
a hyperplane section. The resulting statements can be formulated for higher homotopy groups 
as well, but they are easier in that case, and, possibly, less useful. So we will stick to the 
discussion of fundamental groups and just prove Theorem~\ref{t2}.

Let $\Gamma$ be the fundamental group of a Sasaki manifold $M$ of dimension $2n+1\geq 7$. 
We may assume that the Sasaki structure is quasi-regular, and obtain the associated 
orbifold fibration $\pi\colon M\longrightarrow X$ with $X$ a projective orbifold of complex dimension $n$. 
Now let $Y\subset X$ be a generic hyperplane section. We denote the inclusion of $Y$ in $X$ by $i$,
and we let $N\subset M$ be the preimage of $Y$ under $\pi$, i.e.~the total space of the orbifold
circle bundle restricted from $X$ to $Y$. 

Since $X$ arises as the quotient of the quasi-regular Sasaki structure on $M$, the local uniformizing 
groups of $X$ are cyclic and inject in $S^1$. Restricting to the suborbifold $Y\subset X$, we 
have the same local uniformizing groups with the same injections into $S^1$. In other 
words, the restriction of the polarisation of $X$ to $Y$ is again a polarisation in the sense of 
Ross and Thomas~\cite[Definition~2.7]{RT}. Therefore Theorem~\ref{conversestructure} applies, and $N$ is
a smooth Sasaki manifold, rather than just an orbifold.

This situation gives rise to the following commutative diagram of orbifold homotopy exact sequences:
$$
\begin{tikzcd} [column sep={3em,between origins},
  row sep={3em,between origins},]
 \pi_2^{orb}(Y) \arrow[d, "i_*"] & \stackrel{\partial}{\lra} & \pi_1(S^1)\arrow[d, phantom , sloped , "="] & \lra & \pi_1(N) \arrow[d, "i_*"] & \lra & \pi_1^{orb}(Y) \arrow[d, "i_*"]& \lra & 1  \ \\
 \pi_2^{orb}(X) & \stackrel{\partial}{\lra} & \pi_1(S^1) & \lra & \pi_1(M) & \lra & \pi_1^{orb}(X)   & \lra & 1 \ .
\end{tikzcd}
$$
By the Lefschetz hyperplane theorem for orbifolds or Deligne--Mumford stacks~\cite{Halp} the assumption $n\geq 3$
implies that the vertical inclusion-induced map 
$$
i_*\colon \pi_i^{orb}(Y)\lra\pi_i^{orb}(X)
$$
is an isomorphism for $i=1$ and a surjection for $i=2$. Therefore, this diagram shows, by the usual diagram
chase, that $i_*\colon\pi_1(N)\lra\pi_1(M)$ is also an isomorphism. Now $N$ is a Sasaki manifold of real 
dimension two less than the dimension of $M$ having the same fundamental group. This completes the proof of Theorem~\ref{t2}.

\subsection{The Boothby--Wang construction with control on the fundamental group}

We now prove Theorem~\ref{t3}. Let $\Gamma$ be a projective group. 
Taking products with $\CP^1$ in one direction, and using the Lefschetz hyperplane theorem in the other,
one sees that $\Gamma=\pi_1(X)$ where $X$ can be taken as a smooth complex projective variety of (real) dimension $2n$ for any $n\geq 2$. 
Then $X$ is equipped with an integral K\"ahler class $[\omega]$ so that the corresponding line bundle is ample.
The blow-up of $X$ at a point, topologically $X\#\overline{\CP^n}$, can be endowed with the integral K\"ahler class $k[\omega]-E$ where $E$ is the Poincar\'e dual of the exceptional divisor $D$ and $k\in\N$ is large enough, so that the corresponding line bundle is ample on the blowup.

Consider now the Boothby--Wang fibration $M$ over $X\#\overline{\CP^n}$ with Euler class $e=k[\omega]-E$ and the associated long exact sequence
$$
\cdots\longrightarrow\pi_2 (X\#\overline{\CP^n})\stackrel{\partial}{\longrightarrow}\pi_1(S^1)\longrightarrow\pi_1(M)\longrightarrow\Gamma\longrightarrow 0 \ .
$$
The exceptional divisor $D\cong \CP^{n-1}$ contributes a non-torsion spherical class to $H_2(X\#\overline{\CP^n};\Z)$ on which the Euler class evaluates as $\pm 1$. Thus the map $\partial$ is surjective  by~\eqref{EqMapFact}, that is, we get the desired isomorphism 
$\pi_1(M)\cong \Gamma$.
This completes the proof of Theorem~\ref{t3}.

\begin{remark}\label{orbifoldremark}
This proof generalizes to the orbifold setting in the following way. Instead of starting with a smooth projective 
variety, we may start with any cyclic projective orbifold with an integral polarisation.
We first blow up a smooth
point and modify the polarisation to one on the blowup which evaluates as $\pm 1$ on a projective line in the 
exceptional divisor. Then we apply Theorem~\ref{conversestructure} to this polarisation on the blowup to 
obtain a Sasaki manifold whose fundamental group is the orbifold fundamental group of the cyclic orbifold 
we started with.
\end{remark}

Another variation on the proof of Theorem~\ref{t3} is given by the following:
\begin{prop}
Every finitely presentable group is the fundamental group of a closed $K$-contact manifold of dimension $2n+1$
for any $n\geq 2$.
\end{prop}
Without the control on the dimension, this result appears in~\cite[Theorem~5.2]{BFMT}. The proof
given there is different from ours, and does not prove the case $n=2$. The following argument implements a 
suggestion of~\cite[Remark~2.2]{AK}.
\begin{proof}
Fix an $n\geq 2$ and a finitely presentable group $\Gamma$. 
By a celebrated theorem of Gompf~\cite{Gompf} there exists a closed symplectic 
$2n$-manifold $Y$ such that $\pi_1(Y)=\Gamma$. Since non-degeneracy is an open condition in the space of closed forms,
there exists a symplectic form on $Y$ representing a rational class in cohomology. After multiplication with a large integer we may assume that $Y$ is equipped with a symplectic form $\omega$ representing a primitive integral class $[\omega]$.

After possibly replacing $Y$ by its symplectic blow-up at a point, there exists a spherical class in $H_2(Y;\Z)$ on which $[\omega]$ evaluates as $\pm 1$. Therefore, the Boothby--Wang fibration over $Y$ with Euler class $[\omega]$ is a compact $(2n+1)$-dimensional 
K-contact manifold $N$ with $\pi_1(N)=\Gamma$. 
\end{proof}

\section{Some non-K\"ahler Sasaki groups}\label{s:SasakiNotKaehler}

In this section we construct interesting Sasaki groups which are not K\"ahler, and use some of 
them to prove Theorem~\ref{t1}.

Since $\Sg_5$ is the richest class of Sasaki groups, it makes sense to start in dimension $5$.

\begin{prop}\label{tAlb}
Let $X$ be an aspherical algebraic surface, and $\pi\colon M\longrightarrow X$ the Sasaki 
$5$-manifold obtained by the Boothby--Wang construction using as Euler class $[\omega ]$ the first
Chern class of an ample line bundle on $X$.
If $X$ has complex Albanese dimension $2$, then $\pi_1(M)$ is not a K\"ahler group.
\end{prop}
\begin{proof}
Suppose $Y$ is a compact K\"ahler manifold with fundamental group $\pi_1(Y)=\pi_1(M)$.
Consider the composition
$$
Y \stackrel{c_Y}{\longrightarrow} B\pi_1(Y) = M \stackrel{\pi}{\longrightarrow} X \ ,
$$
where $c_Y$ is the classifying map of the universal covering of $Y$. Both maps 
induce isomorphisms on $H^1 (-;\R)$. However, any four-fold cup product of classes
in $H^1 (X;\R)$ lands in $H^4(X;\R)$, which is spanned by $[\omega ]^2$.
Since $[\omega ]$ is killed by $\pi^*$, so is $[\omega ]^2$. Therefore, the cup-length 
of degree one classes on $Y$ is strictly less than $4$. This means that the 
image of $Y$ under its Albanese map $\alpha_Y\colon Y\longrightarrow Alb(Y)$ 
is a curve $D$, cf.~\cite[p.~23]{ABCKT}.  
A standard argument then implies that $D$ is smooth, and the 
Albanese map has connected fibres, compare~\cite[p.~289]{F}. 
The Albanese image $D$ is necessarily of genus $\geq 2$,
since the first Betti number of $X$, and therefore of $Y$, is at least $4$. 

Since the Albanese map has connected fibres, it induces a surjective homomorphism 
$$
(\alpha_Y)_*\colon\pi_1(Y)\longrightarrow\pi_1(D) \ .
$$
However, as $D$ is of genus $g\geq 2$, its fundamental group has 
trivial centre, and so this homomorphism factors through the quotient map 
$\pi_*\colon\pi_1(Y)=\pi_1(M)\longrightarrow\pi_1(X)$ as follows:

\begin{equation}\label{AlbDim}
\xymatrix{1 \ar[r] & \Z \ar[r] & \pi_1(Y)\ar[d]_{(\alpha_Y)_*}  \ar[r]^{\pi_*} &  \pi_1(X) \ar[dl]  \ar[r] & 1 \ .\\ 
 & & \pi_1(D) 
}
\end{equation}
Since both $(\alpha_Y)_*$ and $\pi_*$ induce isomorphisms on $H^1 (-;\R)$, so does the 
diagonal map on the bottom right. This means that $H^1 (X;\R)$ comes from the curve $D$,
and so the cuplength is $2$, and not $4$. This contradicts the assumption that $X$ has 
real Albanese dimension $4$, and this contradiction finally shows that $Y$ does not exist,
and $\pi_1(M)$ is not a K\"ahler group.
\end{proof}

There are many examples where this result applies.
\begin{ex}
If we take $X$ to be an Abelian surface, we can arrange $\pi_1(M)$ to be the five-dimensional
integral Heisenberg group. That this is not a K\"ahler group was first proved by
Carlson and Toledo~\cite{CT2} using a completely different argument.
\end{ex}

\begin{remark}
Carlson and Toledo~\cite{CT2} pointed out that the integral Heisenberg groups 
of dimension at least $5$ are known to be $1$-formal, so there is no obstruction
there to them being K\"ahler. Indeed, even if this had not been known thirty years 
ago, it now follows via Kasuya's result~\cite[Theorem~1.1]{Kas1} from the fact 
that they are all Sasaki groups in dimension $\geq 5$. The claim in~\cite[Section~3.1]{BFMT}
that these groups are not $1$-formal is wrong.
\end{remark}

\begin{ex}
We can also take for $X$ any Cartesian product of curves of positive genus,
not necessarily of genus one. Then $\pi_1(M)$ is no longer nilpotent as in 
the previous example.
\end{ex}

\begin{ex}
Proposition~\ref{tAlb} applies whenever $X$ is a Kodaira fibration whose Albanese
image is not a curve. Almost all examples of Kodaira fibrations constructed explicitly
do have this property, compare the discussion in~\cite{Bregman} and the references given there.
\end{ex}

Finally, we can also apply Proposition~\ref{tAlb} to ball quotients of complex dimension $2$
whose Albanese image is also of dimension $2$. However, in this case we can
prove more, in that we can allow ball quotients of arbitrary dimension, and we 
can dispense with the assumption about the Albanese dimension.

Let $G\subset PU(n,1)$ be a torsion-free cocompact lattice in the isometry group of the 
$n$-dimensional complex ball $\C H^n$ equipped with the Bergman metric. Then 
$X=\C H^n/G$ is a smooth projective variety, since the K\"ahler class $[\omega ]$ 
of the Bergman metric is the first Chern class of an ample line bundle. 
Let $\pi\colon M\longrightarrow X$ be the principal circle bundle with Euler class $[\omega ]$.
Its total space carries a Sasaki structure obtained by the Boothby--Wang construction. 
Since $X$ is aspherical, the fundamental group $\Gamma = \pi_1(M)$ is the central extension 
\begin{equation}\label{Gamma}
1 \longrightarrow \Z \longrightarrow \Gamma \longrightarrow G \longrightarrow 1 
\end{equation}
with extension class $[\omega ]$.

The following statement is known to some experts, cf.~for example~\cite{Rez}, 
but since a quick and clean proof is not readily available in the literature, we give one here.

\begin{prop}\label{GammanotKaehler}
The extension $\Gamma$ is not a K\"ahler group.
\end{prop}
\begin{proof}
In case that $n=1$, the manifold $M=B\Gamma$ has positive first Betti number,
but the cup product from $H^1$ to $H^2$ vanishes identically. Therefore, $\Gamma$
cannot be K\"ahler by the Hard Lefschetz Theorem.

We now assume $n\geq 2$ and argue by contradiction. So suppose that $Y$ is 
a compact K\"ahler manifold with fundamental group $\Gamma$. Consider the 
classifying map of its universal covering composed with the projection from $M=B\Gamma$
to $X=BG$:
$$
Y \stackrel{c_Y}{\longrightarrow} M  \stackrel{\pi}{\longrightarrow} X \ .
$$
Since $X$ is negatively curved, the composition $\pi\circ c_Y$ is homotopic
to a harmonic map $h\colon Y \longrightarrow X$ by the Eells--Sampson theorem.
Note that on $\pi_1(Y)$ the induced map $h_*$ is just the projection $\Gamma\longrightarrow G$
with kernel $\Z$. The assumption $n\geq 2$ implies that $h_*$ cannot factor through a 
surface group. We now invoke the classification of harmonic maps from compact 
K\"ahler manifolds to ball quotients due to Carlson and Toledo~\cite[Theorem~7.2]{CT}.
This says that either the rank of the differential $Dh$ is everywhere $\leq 2$, in which 
case $h$ has image a closed geodesic or factors through a Riemann surface, or it 
is holomorphic (after perhaps conjugating the complex structure on $X$). The first 
case is not possible because of what we know about $h_*$ at the fundamental group 
level, and so we conclude that $h$ is indeed holomorphic. Since it is non-constant,
its image is a positive-dimensional analytic subvariety of $X$, on which the appropriate 
power of the K\"ahler class $[\omega ]$ is positive. This means $h^*[\omega]\neq 0$,
contradicting the fact that $h^* = c_Y^*\circ\pi^*$, and $[\omega]\in \ker (\pi^* )$
by the Gysin sequence. This contradiction proves the claim.
\end{proof}

The next Proposition completes the proof of Theorem~\ref{t1}.
\begin{prop}
Let $N$ be any compact Sasaki manifold whose fundamental group is the extension
$\Gamma$ from Proposition~\ref{GammanotKaehler}.
Then $\dim (N)\leq 2n+1$.
\end{prop}
\begin{proof}
As usual we may assume the Sasaki structure to be quasi-regular, and consider the 
quotient map $N\longrightarrow Z$, where $Z$ is a projective orbifold. We obtain the 
corresponding exact sequence 
\begin{equation}\label{extC}
1 \longrightarrow C \longrightarrow \Gamma \longrightarrow \pi_1^{orb}(Z) \longrightarrow 1 \ .
\end{equation}
Since $\Gamma$ is torsion-free, $C$ cannot be non-trivial and finite. Moreover, if $C$
were trivial, then $\Gamma= \pi_1^{orb}(Z) $, and using torsion-freeness again, this would
show that $\Gamma$ is projective, contradicting Proposition~\ref{GammanotKaehler}. Thus $C$ must 
be infinite cyclic. Since the centre $C(\Gamma )$ of $\Gamma$ is infinite cyclic, consisting 
of the copy of $\Z$ on the left in~\eqref{Gamma}, we conclude that $C$ is a finite 
index subgroup of $C(\Gamma )$.

Assume first that $C=C(\Gamma )$. Then $\pi_1^{orb}(Z)=G$, and since $G$ is torsion-free,
we conclude $\pi_1(Z)=G$. Moreover, Remark~\ref{euler} shows that the circle bundle 
$N\longrightarrow Z$ is the pullback of $M\longrightarrow X$ under the classifying map
$c_Z\colon Z\longrightarrow X$. On $X$, the $(n+1)^{\textrm{st}}$ power of the Euler class vanishes
for dimension reasons, but this Euler class pulls back to a K\"ahler 
class on $Z$. Therefore, $\dim_{\C} (Z)\leq n$, which implies $\dim (N)\leq 2n+1$.

If $C\subset C(\Gamma )$ has index $k>1$, then we divide~\eqref{Gamma} by $C$,
and obtain 
$$
1 \longrightarrow C(\Gamma)/C=\Z_k \longrightarrow \Gamma/C = \pi_1^{orb}(Z)  \longrightarrow G \longrightarrow 1 \ .
$$
Again the extension class of~\eqref{extC} in $H^2(\pi_1^{orb}(Z);\Z)\subset H^2(Z;\Z)$ is
a pullback from $BG=X$, showing that its $(n+1)^{\textrm{st}}$ power vanishes. Since it is 
a K\"ahler class, we obtain the same dimension bound as before.
\end{proof}

\section{Further restrictions on Sasaki groups}\label{s:further}

In this section we prove Theorems~\ref{t4} and~\ref{t5}, and we discuss some applications and variations. 

\subsection{About Theorem~\ref{t5}}
For the proof we use the following lemma.
\begin{lem}[\cite{KotAIF}]\label{LemFinIndSubgroup}
	Let $\Gamma_1$ and $\Gamma_2$ be two groups. Assume $f_i\colon \Gamma_i\lra Q_i$ is a non-trivial quotient with  kernel $K_i$ and $\vert Q_i\vert=m_i<\infty$ for $i=1,2$. Then the free product $\Gamma_1*\Gamma_2$ admits a finite index subgroup with odd first Betti number.
\end{lem}
\begin{proof}
	Consider the following composition
	$$
	\Gamma_1*\Gamma_2\xlongrightarrow{\pi_{ab}}\Gamma_1\times\Gamma_2\xlongrightarrow{f_1\times f_2}Q_1\times Q_2\ .
	$$
	By the Kurosh subgroup theorem, the kernel of this homomorphism has the form $F_m* K$ where $F_m$ is the free group 
	on $m=(m_1-1)(m_2-1)$ generators and $K$ is a free product of subgroups isomorphic to the $K_i$.
	Now let $f\colon F_m\lra Q$ be a finite quotient with $\vert Q\vert= d$.
	Extend $f$ trivially on $K$ to get a homomorphism $\bar f\colon F_m* K\lra Q$.
	Then the kernel of $\bar f$ has the form $F_n* K*\cdots* K$ where $n=1+d(m-1)$ and $K$ appears $d$ many times.
	Thus, $\ker(\bar f)$ is a finite index subgroup in $\Gamma_1*\Gamma_2$ and 
	$$\betti_1(\ker(\bar f))=n+d\betti_1(K)=1+d(m-1+\betti_1(K))\ .$$
	By picking $d$ even we get a finite index subgroup of $\Gamma_1*\Gamma_2$ with odd first Betti number. 
\end{proof}

We are now ready to prove Theorem~\ref{t5}.
Clearly part~\eqref{johnsonparta} follows directly from Lemma~\ref{LemFinIndSubgroup}, so we only have to prove 
part~\eqref{johnsonpartb}. 
Set $\Gamma=(\Gamma_1*\Gamma_2)\times H$. The proof is divided into two cases.

\subsubsection*{Case 1: $\betti_1(H)$ is even}

By Lemma~\ref{LemFinIndSubgroup}, there exists a finite index subgroup $\Delta\subset\Gamma_1*\Gamma_2$ 
with $\betti_1(\Delta)$ odd.
Hence, the group $\Delta\times H$ is a finite index subgroup of $\Gamma$ with odd first Betti number. 
Thus $\Gamma$ cannot be Sasaki.

\subsubsection*{Case 2: $\betti_1(H)$ is odd}
In this case $\betti_1(\Gamma_1* \Gamma_2)>0$. 
Then we can assume that the first Betti number of $\Gamma_1$ is positive, and so $\Gamma_1$ has finite
quotients of arbitrarily large order.
The first step in the proof of Lemma~\ref{LemFinIndSubgroup} provides a finite index subgroup of $\Gamma_1\times\Gamma_2$ 
of the form $F_m* G$. Moreover, the rank $m$ of $F_m$ can be chosen to be arbitrarily large by using suitable
quotients of $\Gamma_1$ (and some fixed quotient of $\Gamma_2$).
Since the class of Sasaki groups is closed under taking finite index subgroups, we can assume 
$\Gamma=(F_m* G)\times H$ with $m>\betti_1(H)$. 

Let $M$ be a compact Sasaki manifold with $\pi_1(M)=\Gamma$. Consider a quasi-regular Sasaki structure 
$\pi\colon M\lra X $ and let
$$
0\longrightarrow C\longrightarrow\Gamma\longrightarrow\pi_1^{orb}(X )\longrightarrow 0 
$$
be the associated central extension.
Since $C$ is mapped to the centre of $\Gamma$, it must be mapped into $H$. It follows that 
$$
\pi_1^{orb}(X ) = (F_m* G)\times (H/C) \ .
$$

Now $H^1(\pi_1^{orb}(X );\R)$ is endowed with a skew-symmetric non-degenerate bilinear pairing by 
Lemma~\ref{LemNonDegPairingOrbi}. 
Moreover, this factors through the cup product, i.e. 
\begin{equation}\label{EqInducedProductOrbiFundGroup}
	\begin{tikzcd}[row sep=tiny]
		H^1(\pi_1^{orb}(X );\R)\times H^1(\pi_1^{orb}(X );\R) \arrow[r, "\cup"]  &  H^2(\pi_1^{orb}(X );\R)  \arrow[r]  &\R\ .     
	\end{tikzcd} 
\end{equation}
Therefore, for this pairing $H^1(F_m)$ is an isotropic subspace in $H^1(\pi_1^{orb}(X ))\cong H^1(F_m)\oplus H^1(G)\oplus H^1(H/C)$ 
which is orthogonal to $H^1(G)$.
Since $\betti_1(H/C)=\betti_1(H)$, the inequality $m>\betti_1(H)=\betti_1(H/C)$ contradicts the non-degeneracy of the 
skew-symmetric pairing. This contradiction proves that $\Gamma$ is not a Sasaki group. We have now
completed the proof of Theorem~\ref{t5}.

\subsection{About Theorem~\ref{t4}.}
We now go through the different cases in  Theorem~\ref{t4} giving the proofs, and, in some cases, 
further applications.

\subsubsection*{Statement~\ref{nocentre}}

Triviality of the centre means that $C$ has to be trivial. Torsion-freeness then implies that
the kernel $K$ in~\eqref{EqDiagramSasakiGroups} is trivial, and so $\Gamma$ is isomorphic 
to the projective group $\pi_1(X)$. This completes the proof in this case.

Napier and Ramachandran~\cite{Napier} proved that the Thompson group $F$ and its generalisations $F_{n,\infty}$ and $F_n$ 
are not K\"ahler. Since these groups are torsion-free with trivial centre, we obtain:
\begin{cor}
	The Thompson group $F$ and its generalisations $F_{n,\infty}$ and $F_n$ are not Sasaki groups.
\end{cor}

\subsubsection*{Statement~\ref{hyperbolic}}

Since Sasaki groups have even first Betti numbers, an infinite Sasaki group cannot be cyclic.
Therefore we may assume that our hyperbolic group is not 
virtually cyclic. Then by~\cite[Lemma~3.5]{KL} its centre is finite. By torsion--freeness the centre
is trivial, and so Statement~\ref{nocentre} applies.

\subsubsection*{Statement~\ref{Schreier}}

Recall from~\cite[Definition~3.1]{HarpeKot} that a group is Schreier if every finitely generated normal 
subgroup is either finite or of finite index. For a Sasaki group $\Gamma$ 
being Schreier means that $C$ cannot be infinite cyclic. For it it were so, then it would have to be of 
finite index, so that $\Gamma$ would be virtually infinite cyclic. Once $C$ cannot be infinite, torsion-freeness 
again proves that  $\Gamma$ is isomorphic to the projective group $\pi_1(X)$.

Llosa~Isenrich~\cite[Corollary~3.2.9]{Llosa} proved that any torsion-free Schreier and K\"ahler group
with virtually positive first Betti number is an orientable surface group of genus $\geq 2$. Thus we
obtain:
\begin{cor}
Any torsion-free Sasaki and Schreier group with virtually positive first Betti number is 
an orientable surface group of genus $\geq 2$. 
	\end{cor}

\subsubsection*{Statement~\ref{rankone}}

Let $N$ be a closed Riemannian manifold of non-positive curvature with
fundamental group $\pi_1(N)=\Gamma$ a Sasaki group.
By the Cartan--Hadamard theorem, $\Gamma$ is torsion-free. Since $\Gamma$ is Sasaki it
cannot be cyclic, so we may assume that $N$ is of dimension $\geq 2$.

The notion of rank for an abstract group was introduced by Ballmann and Eberlein~\cite{BE},
who proved that for the fundamental groups of closed manifolds of non-positive
sectional curvature the rank agrees with the geometric rank of the
Riemannian metric defined via spaces of parallel Jacobi fields.
The rank is additive under Cartesian products of manifolds
  respectively direct products of groups, and is invariant under passage to
  finite coverings respectively to finite index subgroups. The
  assumption that $N$ be of rank one therefore implies that $N$ is
  locally irreducible and that $\Gamma$ is irreducible (any direct
  factor of~$\Gamma$ would be infinite because the group is
  torsion-free).

  The irreducibility of~$N$ and the assumption~$\dim N\geq 2$ imply
  that $N$ has no Euclidean local de Rham factor, so that by the
  result of Eberlein~\cite[p.~210f]{EbJDG}, the centre of~$\Gamma$
  is trivial. Therefore Statement~\ref{nocentre} applies to give the proof
  in this case.

Of course, if the sectional curvature of $N$ is strictly negative, then the triviality of the centre
follows directly from Preissmann's theorem. However, there are many more manifolds
of non-positive curvature and rank one than there are negatively curved manifolds.

Since the fundamental groups of closed real hyperbolic manifolds of dimension
$\geq 3$ are known not to be K\"ahler by a result of Carlson and Toledo~\cite[Theorem~8.1]{CT},
we conclude:
\begin{cor}\label{c:hyperbolic}
No fundamental group of a closed real hyperbolic manifold of
dimension $\geq 3$ is Sasaki.
\end{cor}
With the added assumption or arithmeticity, this also appears in the statement of~\cite[Proposition~5.4]{BFMT}.
However, the proof given there, while using some of the same arguments we use,
seems elliptical and not quite to the point.

\subsubsection*{Statement~\ref{SSNT}}

 For $N$ to be a locally symmetric space of non-compact type means
  that its universal covering $\tilde N$ is a globally symmetric space
  without compact or Euclidean factors in its de Rham decomposition.
  Thus $\pi_1(N)$ is torsion-free. The absence of a Euclidean de Rham factor 
  again implies triviality of the centre via  the
  result of Eberlein~\cite[p.~210f]{EbJDG}.
 Therefore Statement~\ref{nocentre}  applies.

\subsubsection*{Statement~\ref{3D}}

Suppose that $\Gamma\in\Sg_5$ is the fundamental group of a $3$-manifold $N$.
Since all finite groups are projective by a classical result of Serre, cf.~\cite[p.~6]{ABCKT}, we
only have to consider infinite groups $\Gamma$.

By a result of Jaco~\cite{Jaco1}, finite presentability of $\Gamma$ implies that we may take $N$ to be compact,
possibly with non-empty boundary. Since this does not change the fundamental group, we cap off any
spherical boundary components by balls.

Next, we may assume $N$ to be prime. For if it were not prime, its fundamental group $\Gamma$
would be a non-trivial free product. Since $3$-manifold groups are residually 
finite\footnote{Compare~\cite[p.~60]{aschenbrennerfriedlwilton15}.}, Lemma~\ref{LemFinIndSubgroup} 
would show that $\Gamma$ has virtually odd first Betti number, which is impossible since $\Gamma$ is Sasaki.

Moreover, being Sasaki, $\Gamma$ cannot be virtually cyclic, and so the prime manifold $N$ is irreducible by the 
sphere theorem, cf.~\cite{Mil}. This in particular implies that  $N$ is aspherical and $\Gamma$ is torsion-free.
Thus, if the centre of $\Gamma$ is trivial, the conclusion follows from Statement~\ref{nocentre}.

If the centre of $\Gamma$ is not trivial, then $N$ is Seifert fibered, cf.~\cite[Theorem~2.5.5]{aschenbrennerfriedlwilton15}.
After passing to a suitable finite covering, respectively a finite index subgroup of $\Gamma$, we may then
assume that $N$ is a circle bundle over a compact orientable surface $S$. If $S$ has non-empty boundary,
then the Euler class of the circle bundle is trivial, and the total space $N$ has odd first Betti number,
which is impossible for a Sasaki group $\Gamma$. So $S$ must be closed, and the Euler class of the circle
bundle $N\lra S$ must generate $H^2(S)$. Since $\Gamma$ is not virtually cyclic, $S$ has positive 
first Betti number $b_1(S)=b_1(N)$. However, the cup product $H^1(N)\times H^1(N)\lra H^2(N)$ is 
identically zero, and there are non-trivial triple Massey products of classes in $H^1(N)$.
Thus $N$ is not $1$-formal in this case, and so the result of Kasuya~\cite[Theorem~1.1]{Kas1} implies that
$\Gamma$ can actually not be in $\Sg_5$. 

This completes the proof of Statement~\ref{3D}. Combining this statement with the result of~\cite{Kot3}
yields the following conclusion:
\begin{cor}
Any infinite group in $\Sg_5$ that is also the fundamental group of a $3$-manifold is an oriented surface group.
\end{cor}

\end{document}